\documentclass[a4paper]{aomart}

\fancypagestyle{firstpage}
{
  \fancyhf{}
  \cfoot{\footnotesize\thepage}
}

\usepackage{amsmath}
\usepackage{amssymb}
\usepackage[low-sup]{subdepth}

\def \buffer {\vphantom{$0_{abcdefghijklmnopqrstuvwxyz}$}}

\def \A {\mathcal{A}}
\def \N {\mathbb{N}}
\def \R {\mathbb{R}}

\def \cuddle {\hspace{0.1ex}}
\def \* {\vphantom{()}}

\copyrightnote{}

\title
[On the monotonicity of additive representation functions]
{On the monotonicity of additive representation functions}
\author{Pascal Stumpf}

\begin{document}

\raggedbottom

\maketitle

\begin{abstract}
  We study the monotonicity behavior of three slightly differently defined 
  additive representation functions (as initiated by Erd\H{o}s, S\'{a}rk\"{o}zy and S\'{o}s), 
  answering one open question and another one partially, and give a slightly 
  simpler proof for a result of Chen and Tang. 
\end{abstract}

\section{Introduction}

For a set $\A \subset \N_{0} = \N \cup \{0\}$ of non-negative integers and every $n \in \N_{0}$ 
we define, using the notation $\mathrm{card}(\mathcal{S}) = \sum_{\* \, s \, \in \, \mathcal{S}} 1$ for a finite set $\mathcal{S}$, \buffer 
\begin{equation*}
\begin{split}
  r_{1}(\A, n) & = \mathrm{card}\Big{(}\big{\{}(a, b) \in \A \times \A : a + b = n\big{\}}\Big{)} 
  \\[0.0ex]
  r_{2}(\A, n) & = \mathrm{card}\Big{(}\big{\{}(a, b) \in \A \times \A : a + b = n, \, a \leqslant b\big{\}}\Big{)} 
  \\[0.0ex]
  r_{3}(\A, n) & = \mathrm{card}\Big{(}\big{\{}(a, b) \in \A \times \A : a + b = n, \, a < b\big{\}}\Big{)} 
\end{split}
\end{equation*}
as the additive representation functions $r_{1}$, $r_{2}$ and $r_{3}$ belonging to $\A$, which 
count all solutions of the equation $a + b = n$ inside of $\A$ with slightly more 
restrictions as the index of $r$ increases. 

\vspace{1.0ex}

Our starting point are three results of Erd\H{o}s, S\'{a}rk\"{o}zy and S\'{o}s obtained 
in \cite{4} (and a bit later improved by Balasubramanian in \cite{1}) demonstrating the 
surprising different monotonicity behavior of $r_{1}$, $r_{2}$ and $r_{3}$: 

\vspace{2.0ex}

\noindent
{\bf Theorem 0.} (\cite{4}) 
\itshape
Let $\A$ be an infinite set of positive integers. 
\begin{itemize}
  \item[(1)] $r_{1}(\A, n)$\;can be monotone from a certain point on, only if $\A$ contains 
  all the integers from a certain point on. 
  \vspace{1.0ex}
  \item[(2)] $r_{2}(\A, n)$\;cannot be monotone increasing from a certain point on, when 
  $\lim_{\* \,N \, \to \, \infty} \mathrm{card}(\{1, \dots, N\} \setminus \A) / \log(N) = \infty$. 
  \vspace{1.0ex}
  \item[(3)] There is a set $\A$ such that $\N \setminus \A$ is infinite and $r_{3}(\A, n)$ is monotone 
  increasing for all $n \geqslant 0$. 
\end{itemize}
\upshape

\vspace{2.0ex}

\noindent
Later, in his collection of unsolved problems \cite{5}, S\'{a}rk\"{o}zy asked with respect to 
property (1) of Theorem 0, whether or not one can find an infinite set $\A \subset \N_{0}$ 
such that its upper asymptotic density is less than $1$ and $r_{1}(\A, n)$ is monotone 
increasing for almost all $n$, which we can answer positively. 

\vspace{2.0ex}

\noindent
{\bf Theorem 1.1.} 
\itshape
There does exist an infinite set $\A \subset \N_{0}$ such that its natural 
density is $0$ and $r_{1}(\A, n)$ is monotone increasing almost everywhere: \buffer 
\begin{equation*}
  r_{1}(\A, n) \leqslant r_{1}(\A, n + 1) 
\end{equation*}
holds for almost all $n \in \N$. 
\\[1.0ex]
In addition, we can also find a set $\A \subset \N_{0}$ such that $\N_{0} \setminus \A$ is infinite and $r_{1}(\A, n)$ 
is strictly monotone increasing almost everywhere: \buffer 
\begin{equation*}
  r_{1}(\A, n) < r_{1}(\A, n + 1) 
\end{equation*}
holds for almost all $n \in \N$. 
\upshape

\vspace{2.0ex}

\noindent
Until today it remains unknown, whether or not there exists an infinite set $\A$ 
such that $r_{2}(\A, n)$ is monotone increasing from a certain point on, although 
more and more conditions have been collected under which $r_{2}(\A, n)$ cannot be 
monotone increasing (as in \cite{2} and \cite{3}). On the other hand, in their paper \cite{4} 
Erd\H{o}s, S\'{a}rk\"{o}zy and S\'{o}s noted that perhaps a similar construction of a set $\A$ 
as the one for property (3) in Theorem 0 is also possible for $r_{2}(\A, n)$, however 
we can prove this is not possible in the following sense. 

\vspace{2.0ex}

\noindent
{\bf Theorem 1.2.} 
\itshape
If $\A \subset \N_{0}$ is non-empty and $\N_{0} \setminus \A$ is infinite, then $r_{2}(\A, n)$ 
cannot be monotone increasing for all $n \geqslant 0$. 
\upshape

\vspace{2.0ex}

\noindent
Finally, we give an alternative proof of the following result by Chen and Tang 
\cite{3}, in which we do not need property (1) of Theorem 0 for $r_{1}(\A, n)$ anymore. 

\vspace{2.0ex}

\noindent
{\bf Theorem 1.3.} (\cite{3}) 
\itshape
If $\A \subset \N_{0}$, then $r_{2}(\A, n)$ and $r_{3}(\A, n)$ cannot be strictly 
monotone increasing from a certain point on. 
\upshape

\vspace{2.0ex}

\noindent
We would like to mention that illustrating the pairs $(a, b) \in \N_{0} \times \N_{0}$ as points 
$(a + b, a)$ in the plane, such that the corresponding points of all pairs with the 
same sum $a + b = n$ are on one vertical line, was helpful in finding our proofs, 
where for an integer $c \geqslant 0$ not in $\A$ we then remove all points $(c, x)$ and $(x, c)$ 
with $x \in \N_{0}$ lying on two certain lines. 

\newpage

\section{Proofs} 

Before we start with the proofs of all theorems, let us quickly collect the 
following helpful formulas for $r_{1}$, $r_{2}$ and $r_{3}$ in the special case $\A = \N_{0}$. 

\vspace{2.0ex}

\noindent
{\bf Lemma.} 
\itshape
For $n \in \N_{0}$, we have 
\begin{itemize}
  \item[(1)] $r_{1}(\N_{0}, n) = n + 1$, 
  \vspace{1.0ex}
  \item[(2)] $r_{2}(\N_{0}, n) = \lfloor n / 2 \rfloor + 1$, 
  \vspace{1.0ex}
  \item[(3)] $r_{3}(\N_{0}, n) = \lfloor (n - 1) / 2 \rfloor + 1$, 
\end{itemize}
where $\lfloor x \rfloor$ denotes the largest integer not exceeding $x \in \R$. 
\upshape

\vspace{2.0ex}

\noindent
{\bf Proof.} By definition we have \buffer 
\begin{equation*}
  r_{1}(\N_{0}, n) = \mathrm{card}\Big{(}\big{\{}(a, n - a) : a \in \{0, 1, \dots, n\}\big{\}}\Big{)} = n + 1 \,, 
\end{equation*}
and if $n = 2 m$ ($m \in \N_{0}$) is even, we find \buffer 
\begin{equation*}
\begin{split}
  r_{2}(\N_{0}, n) & = \mathrm{card}\Big{(}\big{\{}(a, n - a) : a \in \N_{0}, \, a \leqslant n / 2 = m\big{\}}\Big{)} = m + 1 \,, 
  \\[0.0ex]
  r_{3}(\N_{0}, n) & = \mathrm{card}\Big{(}\big{\{}(a, n - a) : a \in \N_{0}, \, a < n / 2 = m\big{\}}\Big{)} = m \,, 
\end{split}
\end{equation*}
or when $n = 2 m + 1$ is odd, then \buffer 
\begin{equation*}
\begin{split}
  r_{2}(\N_{0}, n) & = \mathrm{card}\Big{(}\big{\{}(a, n - a) : a \in \N_{0}, \, a \leqslant n / 2 = m + 1 / 2\big{\}}\Big{)} = m + 1 \,, 
  \\[0.0ex]
  r_{3}(\N_{0}, n) & = \mathrm{card}\Big{(}\big{\{}(a, n - a) : a \in \N_{0}, \, a < n / 2 = m + 1 / 2\big{\}}\Big{)} = m + 1 \,, 
\end{split}
\end{equation*}
and both cases together also lead us to the formulas (2) and (3). \hfill $\Box$

\vspace{2.0ex}

\noindent
{\bf Proof of Theorem 1.1.} 
\\[1.0ex]
First, let us choose $\A = \{2^{\cuddle{}i} : i \in \N\}$ whose natural density \buffer 
\begin{equation*}
  \lim_{\* N \, \to \, \infty} \frac{\mathrm{card}(\A \cap \{1, \dots, N\})}{N} = \lim_{\* N \, \to \, \infty} \frac{\lfloor \log_{2}(N) \rfloor}{N} \leqslant \lim_{\* N \, \to \, \infty} \frac{\log_{2}(N)}{N} = 0 
\end{equation*}
does exist and equals $0$. Out of the $\lfloor \log_{2}(N) \rfloor$ members of $\A$ up to $N \geqslant 1$ we 
can build no more than $\log_{2}(N)^{2}$ pairwise sums, or in other words, there are 
at least $N - \log_{2}(N)^{2}$ positive integers $n \leqslant N$ such that \buffer 
\begin{equation*}
  r_{1}(\A, n) = 0 \leqslant r_{1}(\A, n + 1) \,, 
\end{equation*}
and hence the probability that a positive integer $n$ chosen at random satisfies 
$r_{1}(\A, n) \leqslant r_{1}(\A, n + 1)$ is at least \buffer 
\begin{equation*}
  \lim_{\* N \, \to \, \infty} \frac{N - \log_{2}(N)^{2}}{N} = 1 - \lim_{\* N \, \to \, \infty} \frac{\log_{2}(N)^{2}}{N} = 1 - 0 = 1 \,, 
\end{equation*}
as desired. 
\\[1.0ex]
Now, let us choose $\A = \N_{0} \setminus \{2^{\cuddle{}i} : i \in \N\}$ whose natural density is $1 - 0 = 1$, 
and define the family of sets $\A_{\cuddle{}j} = \N_{0} \setminus \{2^{\cuddle{}i} : i \in \N, i \leqslant j\}$, where for $j \in \N$ 
we only have removed the first $j$ powers of $2$. 
\\[1.0ex]
If $n \in \{2^{\cuddle{}j} + 1, 2^{\cuddle{}j} + 2, \dots, 2^{\cuddle{}j \, + \, 1}\}$, we have \buffer 
\begin{equation*}
  r_{1}(\A_{\cuddle{}j \, - \, 1}, n) = r_{1}(\N_{0}, n) - 2 \cdot (\cuddle{}j - 1) = n + 1 - 2 \cdot (\cuddle{}j - 1) 
\end{equation*}
since $a + b \leqslant 2^{\cuddle{}j \, - \, 1} + 2^{\cuddle{}j \, - \, 1} = 2^{\cuddle{}j} < n$ for any $a$ and $b$ in $\{2^{\cuddle{}i} : i \in \N, i \leqslant j - 1\}$. 
Moreover, if $n$ is also not of the form $2^{\cuddle{}j} + 2^{\cuddle{}i}$ with $i \in \N$ ($i \leqslant j$), we even find \buffer 
\begin{equation*}
  r_{1}(\A, n) = r_{1}(\A_{\cuddle{}j}, n) = r_{1}(\A_{\cuddle{}j \, - \, 1}, n) - 2 = n + 1 - 2 \cdot j 
\end{equation*}
(while $r_{1}(\A, n) = n + 1 - 2 \cdot (\cuddle{}j - 1)$ in case of $n = 2^{\cuddle{}j} + 2^{\cuddle{}i}$) together with \buffer 
\begin{equation*}
\begin{split}
  r_{1}(\A, n + 1) & = r_{1}(\A_{\cuddle{}j}, n + 1) 
  \\[0.0ex]
  & \geqslant r_{1}(\N_{0}, n + 1) - 2 \cdot j 
  \\[0.0ex]
  & = (n + 1) + 1 - 2 \cdot j > r_{1}(\A, n) 
\end{split}
\end{equation*}
as long as $n + 1 < 2^{\cuddle{}j \, + \, 1}$. Since there are no more than $j$ numbers of the form 
$2^{\cuddle{}j} + 2^{\cuddle{}i}$ from $2^{\cuddle{}j} + 1$ up to $2^{\cuddle{}j \, + \, 1} - 2$, we have found at least $2^{\cuddle{}j} - 2 - j$ numbers 
$n$ in $\{2^{\cuddle{}j} + 1, 2^{\cuddle{}j} + 2, \dots, 2^{\cuddle{}j \, + \, 1}\}$ such that $r_{1}(\A, n) < r_{1}(\A, n + 1)$. 
\\[1.0ex]
In view of the partition \buffer 
\begin{equation*}
  \N = \{1, 2\} \cup \bigcup_{\* j \, = \, 1}^{\* \infty} \{2^{\cuddle{}j} + 1, 2^{\cuddle{}j} + 2, \dots, 2^{\cuddle{}j \, + \, 1}\} \,, 
\end{equation*}
up to an integer $N \geqslant 1$ we then find at most \buffer 
\begin{equation*}
\begin{split}
  & 2 + \sum_{\* j \, = \, 1}^{\* \lfloor \log_{2}(N) \rfloor \, + \, 1} 2 + j 
  \\[0.0ex]
  & \quad \leqslant 2 + (\lfloor \log_{2}(N) \rfloor + 1) \cdot (2 + \lfloor \log_{2}(N) \rfloor + 1) \leqslant 2 + (\log_{2}(N) + 3)^{2} 
\end{split}
\end{equation*}
numbers $n$ such that $r_{1}(\A, n) \geqslant r_{1}(\A, n + 1)$, and so this time the probability 
that a positive integer $n$ chosen at random satisfies $r_{1}(\A, n) < r_{1}(\A, n + 1)$ is 
again at least \buffer 
\begin{equation*}
\begin{split}
  & \lim_{\* N \, \to \, \infty} \frac{N - (2 + (\log_{2}(N) + 3)^{2})}{N} 
  \\[0.0ex]
  & \quad = 1 - \lim_{\* N \, \to \, \infty} \frac{\log_{2}(N)^{2}}{N} + \frac{6 \log_{2}(N)}{N} + \frac{11}{N} = 1 - (0 + 0 + 0) = 1 \,, 
\end{split}
\end{equation*}
as desired. \hfill $\Box$ 

\vspace{2.0ex}

\noindent
{\bf Proof of Theorem 1.2.} 
\\[1.0ex]
Let $c_{1}, c_{2}, c_{3}, \dots\-$ denote the elements of $\N_{0} \setminus \A$ in increasing order. 
\\[1.0ex]
If $c_{1} = 2 m + 1$ $(m \in \N_{0})$ is odd, then \buffer 
\begin{equation*}
\begin{split}
  r_{2}(\A, c_{1} - 1) & = r_{2}(\A \cap \{0, 1, \dots, c_{1} - 1\}, c_{1} - 1) 
  \\[0.0ex]
  & = r_{2}(\{0, 1, \dots, c_{1} - 1\}, c_{1} - 1) 
  \\[0.0ex]
  & = r_{2}(\N_{0}, c_{1} - 1) = \lfloor (2 m + 1 - 1) / 2 \rfloor + 1 = m + 1 \,, 
\end{split}
\end{equation*}
while on the other side \buffer 
\begin{equation*}
\begin{split}
  r_{2}(\A, c_{1}) & = r_{2}(\A \cap \{0, 1, \dots, c_{1} - 1, c_{1}\}, c_{1}) 
  \\[0.0ex]
  & = r_{2}(\{0, 1, \dots, c_{1} - 1\}, c_{1}) = r_{2}(\N_{0} \setminus \{c_{1}\}, c_{1}) 
  \\[0.0ex]
  & = r_{2}(\N_{0}, c_{1}) - 1 = \lfloor (2 m + 1) / 2 \rfloor + 1 - 1 = m \,, 
\end{split}
\end{equation*}
and so there would be a decrease $r_{2}(\A, c_{1} - 1) > r_{2}(\A, c_{1})$ which means $c_{1}$ has 
to be even, and we distinguish two cases for $c_{1}$. 
\\[1.0ex]
{\bf Case} (1): $c_{1} = 2 x$ for $x > 0$. 
\\[1.0ex]
If the next number $c_{2} = 2 m + 1$ $(m \geqslant x)$ missing from $\A$ is odd, then \buffer 
\begin{equation*}
\begin{split}
  r_{2}(\A, c_{2} - 1) & = r_{2}(\N_{0} \setminus \{c_{1}\}, c_{2} - 1) 
  \\[0.0ex]
  & = r_{2}(\N_{0}, c_{2} - 1) - 1 = \lfloor (2 m + 1 - 1) / 2 \rfloor + 1 - 1 = m \,, 
\end{split}
\end{equation*}
while due to $c_{1} + c_{2} \neq c_{2}$ ($c_{1} > 0$) we get \buffer 
\begin{equation*}
\begin{split}
  r_{2}(\A, c_{2}) & = r_{2}(\N_{0} \setminus \{c_{1}, c_{2}\}, c_{2}) 
  \\[0.0ex]
  & = r_{2}(\N_{0}, c_{2}) - 2 = \lfloor (2 m + 1) / 2 \rfloor + 1 - 2 = m - 1 \,, 
\end{split}
\end{equation*}
and so again there would be a decrease $r_{2}(\A, c_{2} - 1) > r_{2}(\A, c_{2})$ which means 
$c_{2}$ has to be even, and we write $c_{2} = 2 y$ ($y > x$). 
\\[1.0ex]
Assume for a moment that $c_{3}$ is larger than $c_{1} + c_{2} + 1$, then \buffer 
\begin{equation*}
\begin{split}
  r_{2}(\A, c_{1} + c_{2}) & = r_{2}(\N_{0} \setminus \{c_{1}, c_{2}\}, c_{1} + c_{2}) 
  \\[0.0ex]
  & = r_{2}(\N_{0}, c_{1} + c_{2}) - 1 
  \\[0.0ex]
  & = \lfloor (2 x + 2 y) / 2 \rfloor + 1 - 1 
  \\[0.0ex]
  & = x + y \,, 
\end{split}
\end{equation*}
while due to $c_{1} + c_{2} \neq c_{1} + c_{2} + 1$ we get \buffer 
\begin{equation*}
\begin{split}
  r_{2}(\A, c_{1} + c_{2} + 1) & = r_{2}(\N_{0} \setminus \{c_{1}, c_{2}\}, c_{1} + c_{2} + 1) 
  \\[0.0ex]
  & = r_{2}(\N_{0}, c_{1} + c_{2} + 1) - 2 
  \\[0.0ex]
  & = \lfloor (2 x + 2 y + 1) / 2 \rfloor + 1 - 2 
  \\[0.0ex]
  & = x + y - 1 \,, 
\end{split}
\end{equation*}
and this decrease $r_{2}(\A, c_{1} + c_{2}) > r_{2}(\A, c_{1} + c_{2} + 1)$ even remains as long as $c_{3} > c_{2} + 1$, because here \buffer 
\begin{equation*}
  (c_{1} + c_{2}) - c_{i} < (c_{1} + c_{2} + 1) - c_{i} \leqslant (c_{1} + c_{2} + 1) - c_{3} < c_{1} 
\end{equation*}
for $i \geqslant 3$, which means $(c_{1} + c_{2}) - c_{i}$ and $(c_{1} + c_{2} + 1) - c_{i}$ are not in $\N_{0} \setminus \A$, 
and when $c_{j}$ is the largest number less than $c_{1} + c_{2} + 2$ missing from $\A$, then \buffer 
\begin{equation*}
\begin{split}
  r_{2}(\A, c_{1} + c_{2}) & = r_{2}(\N_{0} \setminus \{c_{1}, c_{2}, \dots, c_{j}\}, c_{1} + c_{2}) 
  \\[0.0ex]
  & \geqslant r_{2}(\N_{0}, c_{1} + c_{2}) - 1 - (j - 2) 
  \\[0.0ex]
  & = \lfloor (2 x + 2 y) / 2 \rfloor + 1 - 1 - (j - 2) 
  \\[0.0ex]
  & = x + y - j + 2 \,, 
\end{split}
\end{equation*}
while due to $c_{1} + c_{2} \neq c_{1} + c_{2} + 1$ we get \buffer 
\begin{equation*}
\begin{split}
  r_{2}(\A, c_{1} + c_{2} + 1) & = r_{2}(\N_{0} \setminus \{c_{1}, c_{2}, \dots, c_{j}\}, c_{1} + c_{2} + 1) 
  \\[0.0ex]
  & = r_{2}(\N_{0}, c_{1} + c_{2} + 1) - 2 - (j - 2) 
  \\[0.0ex]
  & = \lfloor (2 x + 2 y + 1) / 2 \rfloor + 1 - 2 - (j - 2) 
  \\[0.0ex]
  & = x + y - j + 1 \,, 
\end{split}
\end{equation*}
and in order to avoid this decrease all that remains is the choice $c_{3} = c_{2} + 1$. 
But then we discover the unavoidable decrease from \buffer 
\begin{equation*}
\begin{split}
  r_{2}(\A, c_{2}) & = r_{2}(\N_{0} \setminus \{c_{1}, c_{2}\}, c_{2}) 
  \\[0.0ex]
  & = r_{2}(\N_{0}, c_{2}) - 2 
  \\[0.0ex]
  & = \lfloor 2 y / 2 \rfloor + 1 - 2 = y - 1 
\end{split}
\end{equation*}
(note that $c_{1} + c_{2} \geqslant 2 + c_{2} > c_{2}$) to \buffer 
\begin{equation*}
\begin{split}
  r_{2}(\A, c_{2} + 1) & = r_{2}(\N_{0} \setminus \{c_{1}, c_{2}, c_{3}\}, c_{2} + 1) 
  \\[0.0ex]
  & = r_{2}(\N_{0}, c_{2} + 1) - 3 
  \\[0.0ex]
  & = \lfloor (2 y + 1) / 2 \rfloor + 1 - 3 = y - 2 
\end{split}
\end{equation*}
(note that $c_{2} + c_{3} > c_{1} + c_{3} > c_{1} + c_{2} \geqslant 2 + c_{2} > c_{2} + 1$). 
\\[1.0ex]
{\bf Case} (2): $c_{1} = 0$. 
\\[1.0ex]
In this case, when $m > 0$ denotes the minimum of $\A$, we can write \buffer 
\begin{equation*}
\begin{split}
  r_{2}(\A, 2 m + n) & = r_{2}(\A \cap \{m + 0, m + 1, \dots, m + n\}, 2 m + n) 
  \\[0.0ex]
  & = r_{2}((\A - m) \cap \{0, 1, \dots, n\}, n) 
  \\[0.0ex]
  & = r_{2}(\A - m, n) 
\end{split}
\end{equation*}
for all $n \geqslant 0$, where for the set $\A - m = \{a - m : a \in \A\}$ (in place of $\A$) we 
have already shown in the first case one can find some $n$ such that \buffer 
\begin{equation*}
  r_{2}(\A - m, n) > r_{2}(\A - m, n + 1) \,, 
\end{equation*}
which in turn also leads to \buffer 
\begin{equation*}
  r_{2}(\A, 2 m + n) > r_{2}(\A, 2 m + n + 1) \,, 
\end{equation*}
and so in any case we have found a decrease. \hfill $\Box$ 

\vspace{2.0ex}

\noindent
{\bf Proof of Theorem 1.3.} 
\\[1.0ex]
Let us suppose there exists an integer $N \geqslant 0$ such that $r_{2}(\A, n)$ and $r_{3}(\A, n)$ 
are strictly monotone increasing for $n \geqslant N$. But then from this point onwards 
$r_{2}(\A, n)$ and $r_{3}(\A, n)$ grow by at least one each whenever $n$ increases by one, 
and thus at $n = 2 N + 3$ we find \buffer 
\begin{equation*}
\begin{split}
  r_{2}(\A, 2 N + 3) & \geqslant r_{3}(\A, 2 N + 3) 
  \\[0.0ex]
  & \geqslant r_{3}(\A, N) + 1 \cdot (N + 3) \geqslant N + 3 
\end{split}
\end{equation*}
in contradiction to \buffer 
\begin{equation*}
\begin{split}
  r_{3}(\A, 2 N + 3) & \leqslant r_{2}(\A, 2 N + 3) 
  \\[0.0ex]
  & \leqslant r_{2}(\N_{0}, 2 N + 3) = \lfloor (2 N + 3) / 2 \rfloor + 1 \leqslant N + 5 / 2 
\end{split}
\end{equation*}
by our Lemma from before. \hfill $\Box$


\begin{thebibliography}{0}

\bibitem{1}{\sc R. Balasubramanian}, \\ {\it A note on a result of Erd\H{o}s, S\'{a}rközy and S\'{o}s}. \\ Acta Arith. {\bf 49} (1987), 45\,--\,53. 

\bibitem{2}{\sc Y.-G. Chen, A. S\'{a}rk\"{o}zy, V. T. S\'{o}s, M. Tang}, \\ {\it On the monotonicity properties of additive representation functions}. \\
Bull. Austral. Math. Soc. (1) {\bf 72} (2005), 129\,--\,138. 

\bibitem{3}{\sc Y.-G. Chen, M. Tang}, \\ {\it On the monotonicity properties of additive representation functions, II}. \\ Discrete Math. {\bf 309} (2009), 1368\,--\,1373. 

\bibitem{4}{\sc P. Erd\H{o}s, A. S\'{a}rk\"{o}zy, V. T. S\'{o}s}, \\ {\it Problems and results on additive properties of general sequences, IV}. \\ Lecture Notes in Math. {\bf 1122}, Springer (1986), 85\,--\,104. 

\bibitem{5}{\sc A. S\'{a}rk\"{o}zy}, \\ {\it Unsolved problems in number theory}. \\ Period. Math. Hungar. {\bf 42} (2001), 17\,--\,35. 

\end{thebibliography}
\end{document}